\documentclass{amsart}
\usepackage{amsthm,amsmath,amsmath,amssymb,amsbsy,graphicx}

\theoremstyle{plain}
\newtheorem{theorem}{Theorem}
\newtheorem{corollary}[theorem]{Corollary}
\newtheorem{proposition}[theorem]{Proposition}
\newtheorem{lemma}[theorem]{Lemma}
\newtheorem{conjecture}[theorem]{Conjecture}

\newtheorem{problem}[theorem]{Problem}

\theoremstyle{definition}

\newtheorem{example}[theorem]{Example}

\theoremstyle{remark}
\newtheorem{remark}[theorem]{Remark}

\newcommand{\R}{\mathbb{R}}

\begin{document}

\title[Reflection Arrangements]{A Geometric
Interpretation of the Characteristic Polynomial of Reflection
Arrangements}

\author{Mathias Drton} 
\address{Department of Statistics \\ 
The University of Chicago}

\author{Caroline J. Klivans}
\address{Departments of Mathematics and Computer Science \\ 
The University of Chicago}

\thanks{First author partially supported by NSF grant DMS-0746265 and
  an Alfred P. Sloan Research Fellowship.}



\begin{abstract}
  We consider projections of points onto fundamental chambers of finite
  real reflection groups.  Our main result shows that for groups of type
  $A_n$, $B_n$, and $D_n$, the coefficients of the characteristic
  polynomial of the reflection arrangement are proportional to the
  spherical volumes of the sets of points that are projected onto faces of
  a given dimension.  We also provide strong evidence that the same
  connection holds for the exceptional, and thus all, reflection groups.
  These results naturally extend those of De Concini and Procesi,
  Stembridge, and Denham which establish the relationship for
  $0$-dimensional projections.  This work is also of interest for the field
  of order-restricted statistical inference, where projections of random
  points play an important role.
\end{abstract}

\maketitle

\section{Introduction}
\label{sec:intro}

A classic problem in statistics is the testing of hypotheses that impose
order-restrictions on a parameter vector or, more generally, require a
parameter vector to belong to a polyhedral cone
$\mathcal{C}\subseteq\mathbb{R}^n$; see for instance
\cite{Robertson1988,Silvapulle2005}.  Let $\pi_\mathcal{C}(x)$ be the
orthogonal projection of $x\in\mathbb{R}^n$ onto $\mathcal{C}$.  If the
projection $\pi_\mathcal{C}(x)$ is in the relative interior of a
$k$-dimensional face of $\mathcal{C}$, then we say that
$\pi_\mathcal{C}(x)$ is $k$-dimensional.  The following problem arises when
studying the probability distributions of statistics suitable for testing
membership in the cone $\mathcal{C}$:

\begin{problem}
\label{prob:1}
Which fraction of the unit sphere in $\mathbb{R}^n$, as measured by
surface volume, is occupied by the points $x$ for which the projection
$\pi_\mathcal{C}(x)$ is $k$-dimensional?
\end{problem}

In statistical calculations, the surface volume fractions appear as weights
in mixtures of probability distributions, most commonly, mixtures of
so-called chi-square distributions.  We denote the surface volume fractions
by $\nu_k$ and refer to them as projection volumes.

\begin{example} 
\label{ex:orthant-plane}
Let $\mathcal{C}=[0,\infty)^2$ be the non-negative orthant in $\R^2$.
All points $x$ in the positive orthant $(0,\infty)^2$ lie inside the cone
and thus have a $2$-dimensional projection $\pi_\mathcal{C}(x)$.  All
points in the non-positive orthant $(-\infty,0]^2$, the polar cone, are
projected to the origin, that is, they have $0$-dimensional projection.
All remaining points being $1$-dimensional, the projection volumes are
$\nu_0=\nu_2=1/4$ and $\nu_1=1/2$.   \qed
\end{example}

In small dimensions, Problem~\ref{prob:1} is resolved easily.  If $n=2$, as
in Example~\ref{ex:orthant-plane}, then the projection volume $\nu_2$ is
determined by the angle between the two rays that span the cone, assuming
the cone is pointed.  The fraction $\nu_1$ is always equal to $1/2$ and
$\nu_2=1-\nu_0-\nu_1$.  If $n=3$, then finding the projection volumes turns
into a problem of spherical trigonometry.  For higher dimensions, however,
the problem becomes more difficult, and statistical practice typically
relies on Monte Carlo integration for calculating the projection volumes.
However, one important exception is known.

\begin{theorem}
\label{thm:stirling}
If $\mathcal{C}=\{x\in\mathbb{R}^n\::\:x_1\le x_2\le \dots \le x_n\}$,
then the projection volumes are proportional to the absolute values of the
Stirling numbers of the first kind, that is, the coefficients of the
polynomial $\chi(t)=t(t-1)(t-2)\cdots (t-n+1)$.
\end{theorem}

In this paper we reprove and generalize this result to two other infinite
families of polyhedral cones.  The generalization is based on the fact that
the cone $\mathcal{C}$ of Theorem~\ref{thm:stirling} can be seen as a
fundamental chamber of the reflection arrangement corresponding to the
reflection group $A_{n-1}$.

Reflection (or Coxeter) groups are classical objects, see for example
\cite{Humphreys, ST}, and have received considerable attention recently
from a more combinatorial perspective, see for
example~\cite{BB,Borovik,Fomin}.  A fundamental chamber of a reflection
arrangement in $\R^n$ is the polyhedral cone obtained by taking the closure of
any one of the arrangement's regions.  To any hyperplane arrangement, one
may associate a polynomial defined over its intersection lattice, known as
the characteristic polynomial of the arrangement.  These notions will be
introduced thoroughly in Section~\ref{sec:reflect}.  Our main results,
Theorems \ref{thm:groups-type-ab} and \ref{thm:groups-type-d} (combined
here) solve Problem~\ref{prob:1} for certain families of reflection
arrangements by relating projection volumes to coefficients of
characteristic polynomials.

\begin{theorem}
\label{thm:ABD} 
Let $\mathcal{W}$ be a reflection group of type $A_n$, $B_n$, or $D_n$
and $\chi(t)$ the characteristic polynomial of the associated reflection
arrangement.  For a generic point $x$ in the fundamental chamber
$\mathcal{C}$ of $\mathcal{W}$, the number of group elements $g\in
\mathcal{W}$ with $k$-dimensional projection $\pi_\mathcal{C}(gx)$ is
equal to the absolute value of the coefficient of $t^k$ in $\chi(t)$.
\end{theorem}

This yields the following Corollary, see Lemma~\ref{lem:probs} and
Remark~\ref{rem:gaussian} of Section~\ref{sec:reflect}.

\begin{corollary}
\label{cor:main}
Let $\mathcal{W}$ be a reflection group of type $A_n$, $B_n$, or $D_n$
and $\chi(t)$ the characteristic polynomial of the associated reflection
arrangement.  Then the projection volumes for any fundamental chamber of
$\mathcal{W}$ are proportional to the absolute values of the coefficients
of $\chi(t)$.
\end{corollary}

We conjecture that Theorem~\ref{thm:ABD} and thus also
Corollary~\ref{cor:main} hold for any finite reflection group.

\begin{conjecture}
\label{conj:main-lemma}
Let $\mathcal{W}$ be a finite reflection group and $\chi(t)$ the
characteristic polynomial of the associated reflection arrangement.  For a
generic point $x$ in the fundamental chamber $\mathcal{C}$ of
$\mathcal{W}$, the number of group elements $g\in \mathcal{W}$ with
$k$-dimensional projection $\pi_\mathcal{C}(gx)$ is equal to the absolute
value of the coefficient of $t^k$ in $\chi(t)$.
\end{conjecture}

As seen in Section~\ref{sec:reflect}, it is sufficient to prove
Conjecture~\ref{conj:main-lemma} for the irreducible reflection groups.  By
Proposition~\ref{prop:rank2}, Conjecture~\ref{conj:main-lemma} is true for
all reflection groups of rank 2.  Theorems~\ref{thm:groups-type-ab} and
\ref{thm:groups-type-d} prove Conjecture~\ref{conj:main-lemma} for the
three infinite families among the irreducible reflection groups of rank at
least 3.  We also offer strong empirical evidence for all but one of
the remaining exceptional irreducible reflection groups.  The one group not
covered by our computer experiments is known as $E_8$, and it is not
covered simply due to its large size.  We strongly believe that
Conjecture~\ref{conj:main-lemma} also holds for $E_8$ and thus all finite
reflection groups.

Recent work by De Concini and Procesi, and Stembridge \cite{DCP} and
subsequent work by Denham \cite{Denham} connects the projection volume for
$0$-dimensional projections to the bottom coefficient of the characteristic
polynomial.  Our Conjecture~\ref{conj:main-lemma} provides a natural
extension.

As we also point out in Remark~\ref{rem:other-interpretations}, the
geometric interpretation of the characteristic polynomial in terms of
projection volumes is different from other interpretations that have
appeared in the literature.  It would thus be interesting to establish
bijections between the different partitions of $\mathcal{W}$ that arise in
the different interpretations.

\subsection*{Acknowledgments}

We are grateful to Bernd Sturmfels for introducing the authors.

\section{Reflection Arrangements}
\label{sec:reflect}

In this section, we define the basic notions from the theory of
hyperplane arrangements and finite reflection groups.  For excellent
references on these topics, we refer the reader to
\cite{stanley:2007} and \cite{Kane}.

\subsection{Hyperplane Arrangements}

A (real central) hyperplane arrangement $\mathcal{A}$ is a collection of
codimension-one linear subspaces of $\mathbb{R}^n$.  All arrangements
appearing in this paper are assumed finite.  The {rank} of an arrangement
$\mathcal{A}$ is defined to be the dimension of the linear space spanned by
the normal vectors to its hyperplanes.  Namely, if $\mathcal{A} = \{H_1,
\ldots, H_m\}$ and $H_i = \{x \in \mathbb{R}^n : \alpha_i x = 0\}$, where
$\alpha_i$ is a non-zero vector in $\mathbb{R}^n$, then
\[
\textrm{Rk}(\mathcal{A}) =  \dim(\textrm{Span}\{\alpha_1, \ldots,
\alpha_m\}).
\]
An arrangement in $\mathbb{R}^n$ is said to be   {essential} if the rank
of $\mathcal{A}$ is $n$.  A   {region} or   {chamber} of
$\mathcal{A}$ is any connected component of the complement of the union of
all the hyperplanes in $\mathcal{A}$.  The closure of a chamber of any
essential arrangement forms a pointed polyhedral cone.  We will consider
projections onto cones formed by particular classes of such hyperplane
arrangements.

Much of the combinatorics of a hyperplane arrangement is encoded by its
intersection lattice.  Given an arrangement $\mathcal{A}$, let
$L(\mathcal{A})$ be the set of all intersections of collections of
hyperplanes in $\mathcal{A}$.  We include $\mathbb{R}^n$ in
$L(\mathcal{A})$ as the intersection of the empty collection.  Define a
partial order on $L(\mathcal{A})$ by reverse inclusion of intersections,
that is, $x \leq y$ in $L(A)$ if $ y \subseteq x$.  Then $L(A)$ forms a
lattice ranked by codimension $n -\dim(x)$.  The M\"{o}bius function $\mu :
L(\mathcal{A}) \rightarrow \mathbb{Z}$ of this lattice is defined
recursively by $\mu(\mathbb{R}^n) = 1$ and $\sum_{z \leq y}\mu(z) = 0$.
The characteristic polynomial of the hyperplane arrangement is the
polynomial
\[
\chi_{A}(t) = \sum_{x \in L(\mathcal{A})} \mu(x) t^{\dim(x)}.
\]
We remark that the Poincar\'{e} polynomial $\pi(\mathcal{A},t)$, as
defined for example in~\cite{OT}, is related to the characteristic
polynomial by
\[
\chi_{\mathcal{A}}(t) =
t^{\textrm{Rk}(\mathcal{A})}\pi(\mathcal{A},-t^{-1}).
\]

\begin{example}
\label{ex:orthant-plane2}
Consider the hyperplane arrangement $\mathcal{A}\subset\mathbb{R}^2$
given by the two coordinate axes $H_1$ and $H_2$.  The intersection
lattice of this arrangement is
$L(\mathcal{A})=\{\mathbb{R}^2,H_1,H_2,\{0\}\}$ with its elements ordered
as $\mathbb{R}^2\le H_1\le\{0\}$ and $\mathbb{R}^2\le H_2\le \{0\}$.  The
M\"obius function thus assigns the values $\mu(\mathbb{R}^2)=1$,
$\mu(H_1)=\mu(H_2)=-1$ and $\mu(\{0\})=1$.  The characteristic polynomial
equals $\chi_\mathcal{A}(t) = t^2-2t+1$.  \qed
\end{example}

\subsection{Reflection Groups}
\label{sec:coxeter}

Let $\mathcal{W}\subset\mathit{GL}(\mathbb{R}^n)$ be a (finite) reflection
group.  Recall that a reflection in $\mathbb{R}^n$ is an isometry which
fixes the points of some hyperplane $H$, often called the mirror of the
reflection.  A finite {reflection arrangement} or {Coxeter arrangement} is
the collection of all mirrors of a finite reflection group.  A fundamental
chamber of the group $\mathcal{W}$ is the closure $\mathcal{C}$ of one of
the regions in the reflection arrangement.  Our main results,
Theorems~\ref{thm:groups-type-ab} and~\ref{thm:groups-type-d} below,
establish a link between the coefficients of the characteristic polynomial
of a reflection arrangement and the projection volumes of its fundamental
chambers.

\begin{example}
\label{ex:orthant3}
The two coordinate axes in the plane discussed in
Example~\ref{ex:orthant-plane2} form a simple example of a reflection
arrangement.  The non-negative cone discussed in
Example~\ref{ex:orthant-plane} is a fundamental chamber for this
arrangement.  The projection volumes $(1/4,1/2,1/4)$ for this chamber are
proportional to the absolute values of the coefficients of the
characteristic polynomial $t^2-2t+1$.  \qed
\end{example}

The following lemma will provide us with a
combinatorial approach to solving Problem~\ref{prob:1} for fundamental
chambers.

\begin{lemma}
\label{lem:probs}
Let $\mathcal{C}$ be a fundamental chamber of a finite reflection group
$\mathcal{W}$.  For $x\in\mathbb{R}^n$, let $b_k(x)$ be the number of
group elements $g\in \mathcal{W}$ for which $\pi_\mathcal{C}(gx)$ is
$k$-dimensional. If $b_k(x)\equiv b_k$ is constant outside a Lebesgue
null set of choices for $x$, then $\pi_\mathcal{C}(X)$ is $k$-dimensional
with probability $b_k/|\mathcal{W}|$ for any random vector $X$ whose
joint distribution is continuous and invariant under the action of
$\mathcal{W}$.
\end{lemma}
\begin{proof}
By the invariance of the distribution of $X$,
\begin{align*}
 P\left(\,\pi_\mathcal{C}(X)\text{ is $k$-dim.}\,\right) &=
 \frac{1}{|\mathcal{W}|}\sum_{g\in \mathcal{W}}
 P\left(\,\pi_\mathcal{C}(gX)\text{ is $k$-dim.}\,\right).
\end{align*}
When summing up conditional probabilities, we obtain that
\[
\sum_{g\in \mathcal{W}} P\left(\, \pi_\mathcal{C}(gX)\text{ is
   $k$-dim.}\,|\,X=x \,\right) = b_k
\]
for almost all $x\in\mathbb{R}^n$. Therefore,
\begin{align*}
 P\left(\,\pi_\mathcal{C}(X)\text{ is $k$-dim.}\,\right) &=
 \frac{1}{|\mathcal{W}|}
 \sum_{g\in \mathcal{W}} \int 
 P\left(\,\pi_\mathcal{C}(gX)\text{ is $k$-dim.}\,|\,
 X=x\,\right)\; dP^X(x)
\end{align*}
is equal to $b_k/|\mathcal{W}|$ as claimed.
\end{proof}

\begin{remark}
\label{rem:gaussian}
A canonical choice for the random vector $X=(X_1,\dots,X_n)$ in
Lemma~\ref{lem:probs} is to pick $X_1,\dots,X_n$ as independent standard
normal random variables.  Then the joint distribution of $X$ is
invariant under the action of the orthogonal group and thus in particular
invariant under the action of all reflection groups in $\mathbb{R}^n$.
Moreover, for this choice, the probability that $X$ is $k$-dimensional is
equal to the projection volume $\nu_k$ that appears in
Problem~\ref{prob:1}.
\end{remark}

\subsection{Irreducible Reflection Groups}

Finite reflections groups coincide with the finite Coxeter groups and as
such there is a classification of the irreducible finite reflection groups;
see for example~\cite{Grove,Humphreys,Kane}.  This classification contains
four infinite families, typically denoted as $A_n$, $B_n=C_n$, $D_n$ and
$I_2(m)$.  As remarked earlier, there has been considerable attention
recently on the combinatorics of Coxeter groups.  A natural combinatorial
perspective of these groups is in terms of symmetries, such a perspective
will be important in our proofs of the main results.  The group $A_n$ is
the symmetry group of the $n$-simplex and thus isomorphic to the symmetric
group on $n+1$ characters.  It acts by permutation of the entries of
vectors in $\mathbb{R}^{n+1}$.  The group $B_n$ is the symmetry group of
the $n$-hypercube and consists of signed permutations of vectors in
$\mathbb{R}^{n}$.  The group $D_n$ is a subgroup of $B_n$ that acts by
signed permutations with an even number of sign changes.  The groups
$I_2(m)$ are the dihedral groups, that is, the symmetry groups of regular
$m$-gons.  In addition to the infinite families, there are $6$ exceptional
cases known as $H_3$, $H_4$, $F_4$, $E_6$, $E_7$ and $E_8$.  In each case
the subscript of the symbol indicates the rank of the associated reflection
arrangement.

Because of the following observation, Conjecture~\ref{conj:main-lemma}
only needs to be proved for irreducible reflection groups.

\begin{lemma}
\label{lem:reduce}
If Conjecture~\ref{conj:main-lemma} is true for the irreducible
reflection groups, then it is true for all reflection groups.
\end{lemma}
\begin{proof}
It suffices to consider two reflection groups
$\mathcal{W}_1\subset\text{GL}(\mathbb{R}^m)$ and
$\mathcal{W}_2\subset\text{GL}(\mathbb{R}^{n-m})$ with fundamental chambers
$\mathcal{C}_1$ and $\mathcal{C}_2$, respectively.  For $j=1,2$, let
$b_{j,k}$ be the number of elements $g\in \mathcal{W}_j$ for which the projection
$\pi_{\mathcal{C}_j}(gx_j)$ of a fixed generic point
$x_j\in\mathcal{C}_j$ is $k$-dimensional.  

Fix a generic point $x$ in the fundamental chamber
$\mathcal{C}=\mathcal{C}_1\times \mathcal{C}_2$ of the reflection group
$\mathcal{W}=\mathcal{W}_1\times
\mathcal{W}_2\subset\text{GL}(\mathbb{R}^n)$.  The $k$-dimensional faces
of $\mathcal{C}$ are of the form $F_1\times F_2$, where $F_1$ is an
$i$-dimensional face of $\mathcal{C}_1$ and $F_2$ is a $(k-i)$-dimensional
face of $\mathcal{C}_2$.  Therefore, the number of elements $g\in
\mathcal{W}$ for which $\pi_\mathcal{C}(x)$ is $k$-dimensional is given by
the convolution
\begin{equation}
 \label{eq:convolution}
 \sum_{i=0}^k b_{1,i}b_{2,k-i}, \quad k=0,\dots,n.
\end{equation}

The reflection arrangement $\mathcal{A}$ associated with $\mathcal{W}$ is
the union of two subarrangements $\mathcal{A}_1$ and $\mathcal{A}_2$
corresponding to $\mathcal{W}_1$ and $\mathcal{W}_2$, respectively.  The
normal vectors for the hyperplanes in $\mathcal{A}_1$ are in
$\mathbb{R}^m\times\{0\}^{n-m}$ and those for $\mathcal{A}_2$ are in
$\{0\}^m\times\mathbb{R}^{n-m}$.  Our claim now follows because it can be
shown that the coefficients of the characteristic polynomial
$\chi_\mathcal{A}(t)$ also obey the convolution rule in
(\ref{eq:convolution}) when expressed in terms of the coefficients of the
characteristic polynomials $\chi_{\mathcal{A}_1}(t)$ and
$\chi_{\mathcal{A}_2}(t)$; compare \cite[Exercise 1.5]{stanley:2007}.
\end{proof}

\begin{example}
Returning to our running example about the non-negative orthant, we see
that the characteristic polynomial factors as $t^2-2t+1=(t-1)^2$.  As
mentioned in the above proof, this factorization mirrors the
factorization of the non-negative orthant as $[0,\infty)\times
[0,\infty)$. \qed
\end{example}

\section{Main Results}
\label{sec:main}

\subsection{Top and Bottom Coefficients}
\label{sec:top-bott-coeff}

The case of the top (or leading) coefficient of the characteristic
polynomial, corresponding to the fundamental chamber $\mathcal{C}$ itself,
is well-known.  The chambers of the arrangement partition space, and the
action of $\mathcal{W}$ on chambers is simply transitive.  Therefore, the
orbit of any generic point in $\mathbb{R}^n$ hits the chamber $\mathcal{C}$
precisely once, and the projection volume for $\mathcal{C}$ is
$1/|\mathcal{W}|$.

The case of the bottom coefficient, the coefficient of the non-zero term of
lowest degree, was shown for reflection groups with so-called
crystallographic root systems by De Concini and Procesi, and
Stembridge~\cite{DCP}.  The general case was proven by
Denham~\cite{Denham}.  The points with lowest-dimensional projection form
the polar cone to $\mathcal{C}$, that is, the cone spanned by the normal
vectors to the reflecting hyperplanes supporting $\mathcal{C}$ (or in the
language of finite root systems, spanned by a set of simple roots).  Hence,
proving the case of the bottom coefficients amounts to determining the
relative spherical volume of the polar cone.  We remark that the work of
Denham~\cite{Denham} involves a combinatorial approach that is case-free in
that it does not use the classification of the irreducible reflection
groups.  His approach rests on the fact that the evaluations of the
characteristic polynomial of a hyperplane arrangement at $\pm 1$ have clear
combinatorial interpretations.

Since the sum of the coefficients of the characteristic polynomial is equal
to the size of the group, we may state the following result.

\begin{proposition}
\label{prop:rank2}
Conjecture~\ref{conj:main-lemma} is true for any reflection arrangement
of rank $2$, and in particular, for the infinite family $I_2(m)$ for
which a canonical choice of a fundamental chamber is
\[
\mathcal{C}(I_2(m)) = \left\{
  x\in\mathbb{R}^2 \::\: x_1\ge 0,\: 
  x_1\cdot\cos\left(\frac{\pi(m-1)}{m}\right)     +
  x_2\cdot\sin\left(\frac{\pi(m-1)}{m}\right)     \ge 0
  \right\}.
\]
\end{proposition}

\subsection{Weighted Projections}
\label{sec:weighted-projections}
In our proof for types $A_n$, $B_n$, and $D_n$, a connection has to be made between projections in
$\mathbb{R}^n$ and $\mathbb{R}^{n-1}$.  To this end, it will be necessary
to work with weighted projections onto the fundamental chamber.  For
$y\in\mathbb{R}^n$, a closed set $\mathcal{K}$ and a weight vector
$\omega\in (0,\infty)^n$, we define the weighted projection
\begin{equation*}
\label{eq:weighted-projection}
\pi_\mathcal{K}(y;\omega) = \arg\min_{x\in\mathcal{K}} \sum_{i=1}^n
\omega_i (x_i-y_i)^2.
\end{equation*}
For example, the weighted projection $\pi_H(y;\omega)$ onto the hyperplane
$H=\{x_i=x_j\}$ is obtained by replacing the $i$-th and $j$-th component of
$y$ with their weighted average
\begin{equation}
\label{eq:weighted-average}
\frac{\omega_i y_i +\omega_{j}y_{j}}{
\omega_i+\omega_{j}}.
\end{equation}

Suppose $r_1,\dots,r_d$ are the normal vectors of a fundamental chamber
$\mathcal{C}$ of a reflection group, that is,
\begin{equation*}
 \label{eq:generic:fundchamber}
\mathcal{C} = \{ x\in\mathbb{R}^n\::\: (r_i,x)\ge 0 \;\text{for all}\;
i\in [d]\}. 
\end{equation*}
Such a set of normal vectors is also known as a set of simple roots,
and they satisfy $(r_i,r_j)\le 0$ for all distinct $i,j\in [d]$; see
\cite[Prop.~4.1.5]{Grove}.  Here, $(x,y)=\sum_i x_iy_i$ denotes the standard
inner product of two vectors $x,y\in\mathbb{R}^n$.  We write
$(x,y)_\omega=\sum_i \omega_ix_iy_i$ for the inner product with respect to
a weight vector $\omega$.

The following Lemma is concerned with points that are on the
``wrong side'' of a supporting hyperplane (or a wall) of
$\mathcal{C}$.  The Lemma shows that when a point is on the ``wrong side''
then its (weighted) projection onto $\mathcal{C}$ can be computed by first
projecting onto the considered wall.

\begin{lemma}
\label{lem:projection-general}
Let $\omega\in (0,\infty)^n$ be a weight vector and
$\phi=(1/\omega_1,\dots,1/\omega_n)$ the vector of inverted weights. If
$(r_i,r_j)_\phi\le 0$ for all $j\not= i$, and $y\in\mathbb{R}^n$
satisfies $(r_i,y)\le 0$, then the projection $\pi_\mathcal{C}(y;\omega)$
is in the wall $H=\{x\in\mathbb{R}^n\::\: (r_i,x)=0\}$ of the fundamental
chamber, and
$\pi_\mathcal{C}(y;\omega)=\pi_\mathcal{C}(\pi_H(y;\omega);\omega)$.
\end{lemma}
\begin{proof}
A well-known fact about projections on closed convex sets states that
$z=\pi_\mathcal{C}(y;\omega)$ if and only if $(y-z,x-z)_\omega\le 0$ for
all $x\in\mathcal{C}$.

Consider a point $z\in\mathcal{C}\setminus H$.  In particular,
$(r_i,z)>0$.  Let $\Omega=\textrm{diag}(\omega)\in\mathbb{R}^{n\times
  n}$.  Pick $\varepsilon>0$ small enough such that
$z_\varepsilon=z-\varepsilon\,\Omega^{-1} r_i$ satisfies
\[
(r_i,z_\varepsilon)=(r_i,z)-\varepsilon\cdot(r_i,r_i)_\phi\ge 0.
\]
Then $z_\varepsilon\in\mathcal{C}$ because
$(r_j,z_\varepsilon)=(r_j,z)-\varepsilon (r_j,r_i)_\phi\ge (r_j,z)\ge 0$
for all $j\not=i$.  Moreover, $(y-z,z_\varepsilon-z)_\omega=
-\varepsilon\cdot( r_i,y-z) >0$ because $(r_i,y)-(r_i,z)<0$.  Therefore,
$\pi_\mathcal{C}(y;\omega)$ lies in the wall $H$.

Since $(r_i,\pi_H(y;\omega))=0$, the previous calculation also implies
that $\pi_\mathcal{C}(\pi_H(y;\omega);\omega)$ is in $H$.  Writing
$\|x\|^2_\omega$ for the norm $(x,x)_\omega$, we have that for any point
$h\in H$,
\begin{equation}
  \label{eq:pythagoras-weighted}
  \|y -h\|^2_\omega = \| y-\pi_H(y;\omega)\|^2_\omega +
  \|\pi_H(y;\omega)-h\|^2_\omega. 
\end{equation}
Since the projection $\pi_\mathcal{C}(y;\omega)$ is in $H$, it can be
determined by minimizing $\| y -h\|^2_\omega$ for $h\in\mathcal{C}\cap
H$.  By (\ref{eq:pythagoras-weighted}), we may instead minimize
$\|\pi_H(y;\omega)-h\|^2_\omega$.  This latter minimization, however,
also yields $\pi_\mathcal{C}(\pi_H(y;\omega);\omega)$ because
$\pi_\mathcal{C}(\pi_H(y;\omega);\omega)\in H$.  Hence,
$\pi_\mathcal{C}(y;\omega)=\pi_\mathcal{C}(\pi_H(y;\omega);\omega)$ as
claimed.
\end{proof}

\begin{example}[Groups of Type $A_n$]
\label{ex:weights-A}
The group $A_{n-1}$ has simple roots $r_i=e_{i+1}-e_i$,
$i=1,\dots,n-1$.  They define the chamber
\[
\mathcal{C}(A_{n-1})=\left\{x\in\mathbb{R}^n\::\: x_1\le \dots \le
  x_n\right\}
\]
from Theorem~\ref{thm:stirling}.  Each root has only two non-zero
entries, and $(r_i,r_{i+1})_\phi=-1/\omega_{i+1}< 0$ and
$(r_i,r_j)_\phi=0$ if $j\ge i+2$.  Hence, the condition of
Lemma~\ref{lem:projection-general} that requires $(r_i,r_j)_\phi\le 0$
for all $i\not= j$ holds for any weight vector $\omega$.  \qed
\end{example}

\begin{example}[Groups of type $B_n$]
\label{ex:weights-B}
The group $B_n$ has simple roots $r_1=e_1$ and $r_i=e_i-e_{i-1}$,
$i=2,\dots,n$.  They define the chamber
\[
\mathcal{C}(B_n)=\left\{x\in\mathbb{R}^n\::\: 0\le x_1\le \dots \le
x_n\right\}. 
\]
Again the condition of Lemma~\ref{lem:projection-general} holds for any
weight vector $\omega$ because $(r_1,r_2)_\phi=-1/\omega_{1}<0$,
$(r_i,r_{i+1})_\phi=-1/\omega_{i}< 0$ and $(r_i,r_j)_\phi=0$ if $j\ge
i+2$.  \qed
\end{example}

\begin{example}[Groups of type $D_n$]
\label{ex:weights-D}
A natural choice of simple roots for the group $D_n$ is $r_1=e_1 +
e_2$, $r_2=e_2 - e_1$ and $r_i=e_i - e_{i-1}$, $i=3,\dots,n$.  They
define the chamber
\[
\mathcal{C}(D_n)=\left\{x\in\mathbb{R}^n\::\: |x_1|\le x_2 \dots \le
  x_n\right\}.
\]
In this case, the requirement that $(r_i,r_j)_\phi\le 0$ for all
$i\not=j$ appearing in Lemma~\ref{lem:projection-general} does present a
condition on a positive weight vector $\omega$.  Namely, it needs to hold
that $\omega_1\ge \omega_2$ because $(r_1,r_2)_\phi = 1/\omega_2 -
1/\omega_1$. \qed
\end{example}

\subsection{Groups of type $A_n$, $B_n$, $D_n$}
\label{sec:main-result}

In this section, we prove Theorem~\ref{thm:ABD}, that is, we show that
Conjecture~\ref{conj:main-lemma} holds for the three infinite families of
real irreducible reflection groups, $A_n$, $B_n$, and $D_n$.
The result is known to hold for the groups $A_n$; recall
Theorem~\ref{thm:stirling}.  This case was first proven in
\cite{miles:1959}.  Our work provides and generalizes a geometric version
of this proof.

\begin{figure}
\centering
\includegraphics[width=3in]{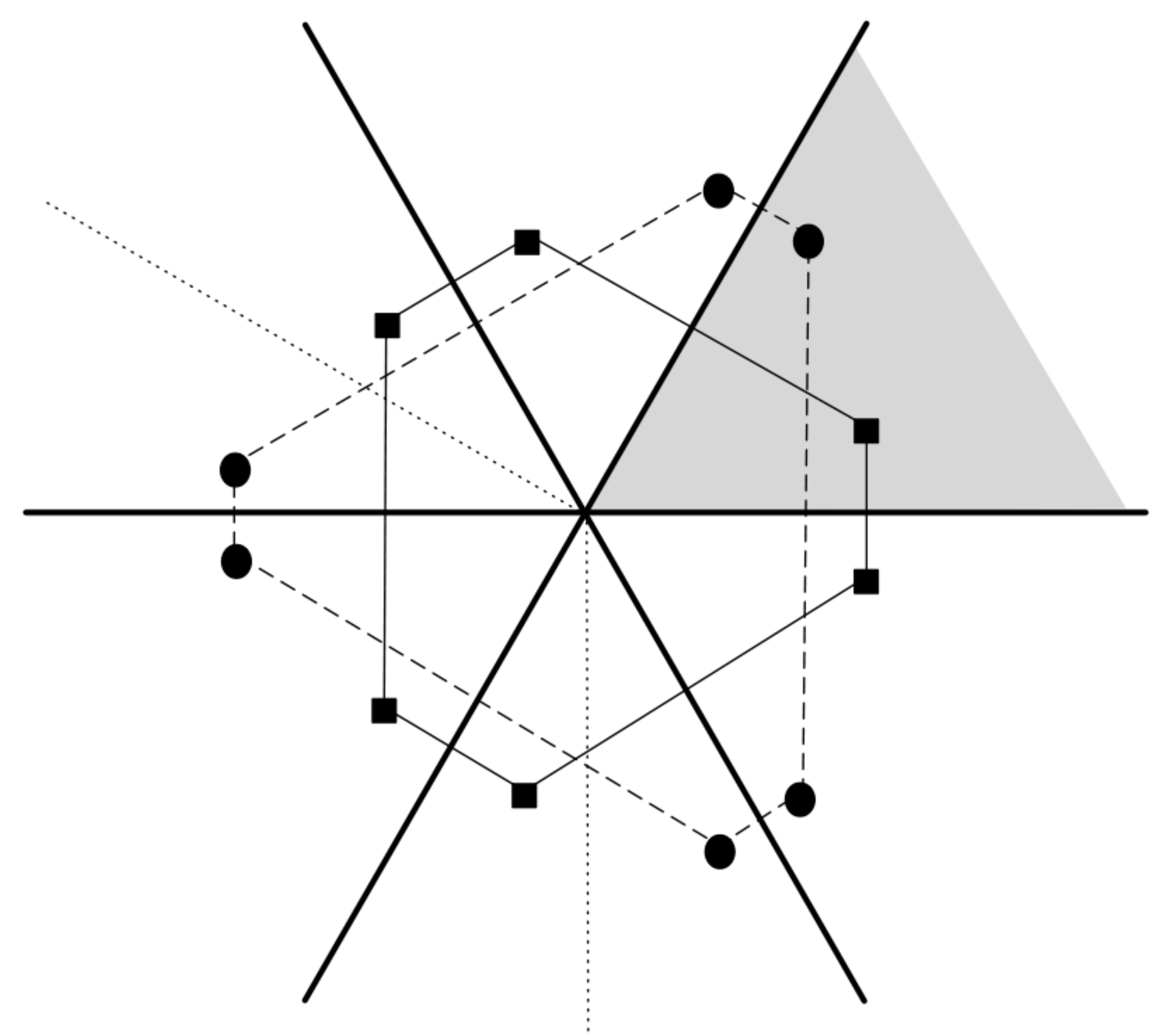}
\caption{Orbits of two points under the reflections of $A_2$.  The
projection chambers are defined by the chosen fundamental chamber
(shaded) and the dotted half lines.  Moving counter-clockwise, the round
point is mapped only once into the first $1$-dimensional projection
chamber, whereas the square point lands in it twice.}
\label{fig:A3}
\end{figure}

One difficulty in proving Theorem~\ref{thm:ABD} is illustrated in
Figure~\ref{fig:A3}, which concerns the group $A_2$.  The figure shows the
orbits of two choices of a generic point $x$.  We see that varying $x$ can
change the dimensionality of the projection of a particular point $gx$ in
the orbit.  Nevertheless, the number of points in the orbit that have a
projection of fixed dimension is independent of the choice of $x$.

Our proof
relies on the well-known factorization of the characteristic polynomial of
a reflection arrangement:
\begin{equation}
\label{eq:factor}
\chi_\mathcal{A}(t) = (t - e_1)(t-e_2) \cdots (t-e_n),
\end{equation}
where the $e_i$s are known as the exponents of the group; see
\cite[Thm.~4.137]{OT}.  For the considered groups, we will be able to build
an induction argument around this factorization by decomposing each group
into $e_n$ pieces that each are in one-to-one correspondence with a
reflection group of lower rank.  Here, $e_n$ is the highest exponent of the
group.

\begin{theorem}
\label{thm:groups-type-ab}
Let $\mathcal{W}$ be a reflection group of type $A_n$ or $B_n$, and
$\omega\in(0,\infty)^n$ be any weight vector.  For a generic point $x$ in the
fundamental chamber $\mathcal{C}$ of $\mathcal{W}$, the number of group
elements $g\in \mathcal{W}$ with $\pi_{\mathcal{C}}(gx;g\omega)$
$k$-dimensional equals the absolute value of the coefficient of $t^k$ in
the characteristic polynomial $\chi_{\mathcal{W}}(t)$.
\end{theorem}
\begin{proof}
We argue by induction on $n$, the index of the group.  The claim is
easily seen to be true for the induction base $n = 1$.  Although the
method of proof for the induction step is the same, we must break up the
two families separately.  

{\em Case $A$}.  The irreducible reflection group $A_{n-1}$ is
isomorphic to the group of permutations of $n$ characters and has the
canonical fundamental chamber $\mathcal{C}(A_{n-1})$ discussed in
Example~\ref{ex:weights-A}.  We abbreviate the projection map onto
this chamber by $\pi_{n-1}$.  The reflection arrangement corresponding to
$A_{n-1}$ is the collection of hyperplanes:
\[
\mathcal{A}_{n-1} = (\,\{x_i - x_j = 0\} \::\: 1 \leq i < j \leq n),
\]
also known as the braid arrangement.  
It has characteristic polynomial
\[
\chi^A_{n-1}(t) = t(t-1)\cdots (t-n+1).
\]

Fix a generic point $x\in\mathcal{C}(A_{n-1})$.  In particular, the
coordinates of $x$ are ordered as $x_1<\cdots< x_n$.    Let
\[
A_{n-1,n} = \left\{ g\in A_{n-1} \::\: (gx)_n = x_n  \right\},
\]
and, for $j=1,\dots,n-1$, 
\[
A_{n-1,j} = \left\{ g\in A_{n-1} \::\: (gx)_i = x_n \;\text{and}\;
  (gx)_{i+1}=x_j \;\text{for some}\; i\in[n]  \right\}.
\]
Clearly, the sets $A_{n-1,1},\dots,A_{n-1,n}$ form a partition of
$A_{n-1}$. In what follows, $\rho_j:\mathbb{R}^n\to\mathbb{R}^{n-1}$
denotes the restriction map omitting the $j$-th coordinate.

The set $A_{n-1,n}$ is a subgroup isomorphic to $A_{n-2}$.  For $g\in
A_{n-1,n}$, define $\bar g$ to be the element of $A_{n-2}$ for which
$\rho_n (gx) = \bar g \rho_n(x)$.  Then $\pi_{n-1}(gx,g\omega)$ is
$k$-dimensional if and only if $\pi_{n-2}(\bar g\rho_n(x),\bar
g\rho_n(\omega))$ is $(k-1)$-dimensional with respect to the chamber
$\mathcal{C}(A_{n-2})$.  By the induction hypothesis, the number of
$k$-dimensional elements in $A_{n-1,n}$ equals the coefficient of $t^k$
in $t\cdot \chi^A_{n-2}(t)$ in absolute value.

If $j<n$, then $A_{n-1,j}$ is again in 1:1-correspondence to $A_{n-2}$.
Let $g\in A_{n-1,j}$ satisfy $(gx)_i = x_n$ and $(gx)_{i+1}=x_j$.  The
corresponding element $\bar g\in A_{n-2}$ is determined by $\rho_i(gx)
= \bar g \rho_n(x)$.  Let $H=\{x_i=x_{i+1}\}$.  Then the $i$-th and the
$(i+1)$-st components in $\pi_H(gx;g\omega)$ are equal; recall
(\ref{eq:weighted-average}).  Since $(gx)_i = x_n> x_j=(gx)_{i+1}$ the
point $gx$ is on the ``wrong side'' of the wall $H$.
Lemma~\ref{lem:projection-general} implies that $\pi_{n-1}(gx,g\omega)$
is $k$-dimensional if and only if $\pi_{n-1}(\pi_H(gx;g\omega),g\omega)$
is $k$-dimensional.  Since $\pi_H(gx;g\omega)=g\pi_{K}(x;\omega)$ for
$K=\{x_j=x_n\}$, this in turn occurs if and only if $\pi_{n-2}(\bar
g\rho_n(\pi_{K}(x;\omega)),\bar g\bar\omega)$ is $k$-dimensional.
The new weights $\bar\omega\in(0,\infty)^{n-1}$ satisfy
$\bar\omega_l=\omega_l$ for $l<j$, $\bar\omega_j=\omega_j+\omega_n$ and
$\bar\omega_l=\omega_{l+1}$ for $j<l<n$.  By the induction hypothesis,
the number of $k$-dimensional elements in $A_{n-1,j}$ equals the absolute
value of the coefficient of $t^k$ in $\chi^A_{n-2}(t)$.

Adding up the count of $k$-dimensional elements in each one of the sets
$A_{n-1,j}$, we obtain the number of $k$-dimensional elements in
$A_{n-1}$.  This number is thus equal in absolute value to the
coefficient of $t^k$ in
\[
t\cdot \chi^A_{n-2}(t) - (n-1)\cdot \chi^A_{n-2}(t) = \chi^A_{n-1}(t).
\]

{\em Case $B$}.  The irreducible reflection group $B_n$ acts by signed
permutations and has the canonical fundamental chamber $\mathcal{C}(B_n)$
discussed in Example~\ref{ex:weights-B}.  We denote the projection map
onto this chamber by $\pi_{n}$. The reflection arrangement corresponding
to $B_n$ consists of the hyperplanes:
\[
\mathcal{B}_n = (\{x_i - x_j = 0\}, \{x_i + x_j = 0\}, \{x_k = 0\} \::\: 1
\leq i < j \leq n, \, \, 1 \leq k \leq n).
\]
The characteristic polynomial for the $B_n$ arrangement is
\[
\chi^B_{n}(t) = (t-1)(t-3) \cdots (t - (2n-1)).
\]

Fix a generic point $x \in \mathcal{C}(B_n)$, which has coordinates
$0<x_1<\dots<x_n$.  For a partition of $B_n$, let
\begin{align*}
  B_{n,n} &= \{g \in B_n \, : \, (gx)_n = x_n \} \textrm{ and}\\
  B_{n,0} &= \{g \in B_n \, : \, (gx)_1 = -x_n \},
\end{align*}
and, for $j = 1, \ldots n-1$, 
\begin{align*}
  B_{n,j} &= 
  \left\{ g \in B_n :   \begin{array}{l} 
      (gx)_i = x_n \textrm{ and } (gx)_{i+1} = x_j,  \textrm{ or }\\
      (gx)_{i+1} = -x_n \textrm{ and } (gx)_i = -x_j
    \end{array}  \textrm{ for some } i \in [n-1]
  \right\}\\
  &\textrm{and}\\
  B_{n,-j} &= 
  \left\{
    g \in B_n  : \begin{array}{l}
      (gx)_i = x_n \textrm{ and } (gx)_{i+1} = -x_j, \textrm{ or }\\
      (gx)_{i+1} = -x_n \textrm{ and } (gx)_i = x_j 
    \end{array} \textrm{ for some } i \in [n-1]
  \right\}.
\end{align*}
The sets $B_{n,0}$, $B_{n,n}$, $B_{n,j}$, $B_{n,-j}$, $1\leq j \leq n-1$,
indeed form a partition of $B_n$.

The set $B_{n,n}$ is a subgroup isomorphic to $B_{n-1}$.  For $g\in
B_{n,n}$, define $\bar g$ to be the element of $B_{n-1}$ for which
$\rho_n(gx) = \bar g \rho_n(x)$.  Then $\pi_{n}(gx,g\omega)$ is
$k$-dimensional if and only if $\pi_{n-1}(\bar g\rho_n(x),\bar
g\rho_n(\omega))$ is $(k-1)$-dimensional with respect to the chamber
$\mathcal{C}(B_{n-1})$.  By the induction hypothesis, the number of
$k$-dimensional elements in $B_{n,n}$ equals the coefficient of $t^k$ in
$t\cdot \chi^B_{n-1}(t)$ in absolute value.

The set $B_{n,0}$ is in 1:1-correspondence with $B_{n-1}$.  Let $H =
\{x_n = 0\}$.  Then $\pi_{H}(x;\omega)$ is a vector with $n$-th coordinate
equal to $0$.  Since $(gx)_1 = -x_n$ for $g \in B_{n,0}$, the first
coordinate of $g\pi_H(x;\omega)$ is $0$.  This follows because $x_n$ is
the largest entry in $x$.  Define $\bar{g}$ to be the element of
$B_{n-1}$ such that $\rho_1(gx) = \bar{g}\rho_n(x)$.
Lemma~\ref{lem:projection-general} implies that $\pi_{n}(gx,g\omega)$ is
$k$-dimensional if and only if $\pi_{n}(g\pi_H(x;\omega),g\omega)$ is
$k$-dimensional.  This in turn occurs if and only if $\pi_{n-1}(\bar
g\rho_n(\pi_H(x;\omega)),\bar g\bar\omega)$ is $k$-dimensional, where the
new weights $\bar\omega\in(0,\infty)^{n-1}$ satisfy
$\bar\omega_l=\omega_{l+1}$ for $l<n$.

For $1 \leq j \leq n-1$, the sets $B_{n,j}$ and $B_{n,-j}$ are again in
1:1-correspondence with $B_{n-1}$.  If $g$ is an element in $B_{n,j}$ or
$B_{n,-j}$, then $g$ satisfies one of the four defining relationships
above.  In each case, the point $x$ is mapped by $g$ to the ``wrong
side'' of the bounding wall $\{x_i = x_{i+1}\}$, namely, $(gx)_i >
(gx)_{i+1}$.  Note also that $(gx)_i$ and $(gx)_{i+1}$ appear with equal
signs in $B_{n,j}$ and opposite signs in $B_{n,-j}$.  Therefore we may
apply Lemma~\ref{lem:projection-general} and average coordinates and add
weights as in the $A$-case. It follows that the number of $k$-dimensional
elements in either $B_{n,j}$ or $B_{n,-j}$ equals the absolute value of
the coefficient of $t^k$ in $\chi^B_{n-1}(t)$.  As there are a total of
$2n -2$ sets of the form $B_{n,j}$ or $B_{n,-j}$, along with the
contributions from $B_{n,n}$ and $B_{n,0}$, we have shown that the total
number of $k$-dimensional elements in $B_n$ is the absolute value of the
coefficient of $t^k$ in
\[
t \cdot\chi^B_{n-1}(t) - (2n-1)\cdot \chi^B_{n-1}(t) = \chi^B_{n}(t),
\]
which was our claim.
\end{proof}

Our second main result concerns the groups of type $D_n$.  As is clear from
Example~\ref{ex:weights-D}, the use of arbitrary weighted projections is no
longer possible.  Nevertheless, we can prove
Conjecture~\ref{conj:main-lemma} in the original unweighted version that is
of most interest.  Compared to the $A$- and the $B$-case, a notable
difference in the proof below is that we do not break up $D_n$ into copies
of $D_{n-1}$.

\begin{theorem}
\label{thm:groups-type-d}
For a generic point $x$ in a fundamental chamber $\mathcal{C}$ of the
reflection group $D_n$, the number of elements $g\in D_n$ that have
$\pi_{\mathcal{C}}(gx)$ $k$-dimensional equals the absolute value of the
coefficient of $t^k$ in the characteristic polynomial $\chi^D_n(t)$.
\end{theorem}
\begin{proof}
The irreducible reflection group $D_n$ acts by signed permutations with
an even number of sign changes.  It has the canonical fundamental chamber
$\mathcal{C}(D_n)$ discussed in Example~\ref{ex:weights-D}.  We denote
the projection onto this chamber by $\pi^D_n$.  The reflection
arrangement corresponding to $D_n$ is the collection of hyperplanes:
\[
\mathcal{D}_n = (\,\{x_i - x_j = 0\}, \{x_i + x_j = 0\} \::\: 1 \leq i < j \leq
n).
\]
The characteristic polynomial of this arrangement is
\[
\chi^D_n(t) = (t-1)(t-3) \ldots (t - 2n+3)(t - n+1).
\]

We now adopt the same partitioning method as in the proof of
Theorem~\ref{thm:groups-type-ab}.  The partition of $D_n$ results from a
coarsening of the decomposition in the $B$-case.  Fix a generic point $x
\in \mathcal{C}(D_n)$.  In particular, $0<|x_1|<x_2<\dots<x_n$.  Let
\[
D_{n,n} = \{g \in D_n \, : \, (gx)_n = x_n \textrm{ or } \, (gx)_1 = -x_n
\}.
\]
For $j = 1, \ldots n-1$, let
\begin{displaymath}
  D_{n,j} = 
  \left\{ 
    g \in D_n  :
    \begin{array}{l} 
      (gx)_i = x_n \textrm{ and } (gx)_{i+1} = x_j, \textrm{ or }\\
      (gx)_{i+1} = -x_n \textrm{ and } (gx)_i = -x_j, \textrm{ or }\\
      (gx)_i = x_n \textrm{ and } (gx)_{i+1} = -x_j, \textrm{ or }\\
      (gx)_{i+1} = -x_n \textrm{ and } (gx)_i = x_j  
    \end{array} 
    \textrm{ for some } i \in [n-1] \right\}.
\end{displaymath}
The sets $D_{n,j}$, $1\le j\le n$, partition $D_n$.

The set $D_{n,n}$ is a subgroup isomorphic to $B_{n-1}$.  The explicit
isomorphism is given by mapping $g\in D_{n,n}$ to $\bar g\in B_{n-1}$
defined as follows.  If $g\in D_{n,n}$ satisfies $(gx)_n=x_n$, then we
define $\bar g\in B_{n-1}$ by requiring that $\rho_n(gx) = \bar g
\rho_n(x)$.  Note that $\bar g$ is a signed permutation with an even
number of sign changes.  If $g\in D_{n,n}$ satisfies $(gx)_1=-x_n$, then
we define $\bar g\in B_{n-1}$ by $\rho_1(gx) = \bar g \rho_n(x)$.  Now
$\bar g$ is a signed permutation with an odd number of sign changes.
With this correspondence between $D_{n,n}$ and $B_{n-1}$ it holds that
$\pi^D_{n}(gx)$ is $k$-dimensional if and only if $\pi^B_{n-1}(\bar
g\rho_n(x))$ is $(k-1)$-dimensional.  Here,
$\pi^B_{n-1}$ stands for the projection onto the fundamental chamber
$\mathcal{C}(B_{n-1})$.  By Theorem~\ref{thm:groups-type-ab}, the
absolute value of the coefficient of $t^k$ in $t\cdot \chi^B_{n-1}(t)$
enumerates the $k$-dimensional elements in $D_{n,n}$.

For $1 \leq j \leq n-1$, the sets $D_{n,j}$ are again in
1:1-correspondence with $B_{n-1}$.  An element $g\in D_{n,j}$ satisfies
one of the four defining relationships above.  In each case, the point
$x$ is mapped by $g$ to the ``wrong side'' of the bounding wall $\{x_i =
x_{i+1}\}$, namely, $(gx)_i > (gx)_{i+1}$.  Note also that $(gx)_i$ and
$(gx)_{i+1}$ may appear with both equal and opposite signs.  Being in an
unweighted situation, or rather in a situation with weight vector
$\omega=(1,\dots,1)$, we may apply Lemma~\ref{lem:projection-general} and
average coordinates and add weights as in the $A$- and the $B$-case.  By
Theorem~\ref{thm:groups-type-ab}, we obtain that the number of
$k$-dimensional elements in $D_{n,j}$ equals the absolute value of the
coefficient of $t^k$ in $\chi^B_{n-1}(t)$.  There being a total of $n -1$
sets of the form $D_{n,j}$, along with the contribution from $D_{n,n}$,
the total number of $k$-dimensional elements in $D_n$ is seen to be the
absolute value of the coefficient of $t^k$ in
\[
t \cdot\chi^B_{n-1}(t) - (n-1)\cdot \chi^B_{n-1}(t) = \chi^D_{n}(t),
\]
as claimed.
\end{proof}

\begin{remark}
 \label{rem:breakdown}
Our method of proof breaks up a Coxeter group into copies of a suitable
subgroup.  One might hope to extend this to a general method for all
reflection groups.  However, this is not possible as can be seen for
example from the fact that removing the highest factor $(t-e_n)$ in the
factorization in (\ref{eq:factor}) need not yield the characteristic
polynomial of a reflection group \cite[Section 6.5]{OT}.  In this context, see
also \cite{Barcelo} where it is shown that if $\mathcal{W}'$ is a
parabolic subgroup of a reflection group $\mathcal{W}$, then the
characteristic polynomial $\chi_{\mathcal{W}'}(t)$ divides $\chi_W(t)$ if
and only if $\mathcal{W}$ is of type $A_n$ or $B_n$, or $\mathcal{W}'$ is
of rank $1$.
\end{remark}

\section{Exceptional Groups}

The proof method used in the previous section decomposes a reflection group
into copies of a reflection subgroup and proceeds inductively.  The method,
however, breaks down for the exceptional irreducible reflection groups;
recall Remark~\ref{rem:breakdown}.  Nevertheless, we believe that
Conjecture~\ref{conj:main-lemma} also holds for these groups.  In the
remainder of this section we describe simulation evidence that supports
this belief.

We tested our conjecture for the groups $H_3$, $H_4$, $F_4$, $E_6$, and
$E_7$.  The only other exceptional irreducible group, $E_8$, was too large
for our implementation.  For each considered group, we randomly chose
points in the fundamental chamber and counted the number of points with
$k$-dimensional projections in the resulting orbits.  We considered 1000
randomly chosen points for the groups of rank at most 6, and 50 points for
$E_7$. Conjecture~\ref{conj:main-lemma} held in all tested cases, that is,
in each test orbit the number of group elements with $k$-dimensional
projection was equal to the absolute values of the coefficient of $t^k$ in
the characteristic polynomial $\chi(t)$.

Suppose the considered fundamental chamber 
\begin{equation*}
 \label{eq:generic:fundchamber:2ndtime}
 \mathcal{C} = \{ x\in\mathbb{R}^n\::\: (r_i,x)\ge 0 \;\text{for all}\;
 i\in [d]\}. 
\end{equation*}
is defined by the simple roots $r_1,\dots,r_d$.  Let $V$ be the linear
space spanned by $\{r_1, \ldots, r_d\}$.  Define the dual roots
$s_1,\dots,s_d\in V$ by requiring that $(s_i,r_j)=\delta_{ij}$.  Then
$\mathcal{C}\cap V$ is the simplicial cone spanned by $\{s_1,\dots,s_d\}$;
see \cite[\S4.2]{Grove}.  Our procedure for generating random points in the
fundamental chamber is based on this spanning set representation, namely,
it picks points uniformly at random from the convex hull of
$\{s_1,\dots,s_d\}$.

Given a point $x\in\mathcal{C}$, we need to visit each point in the orbit
of $x$.  We accomplished this using the {\tt traverse} function from John
Stembridge's Maple package {\tt coxeter} \cite{Stembridge}.  In order to
calculate the dimensions of the projections $\pi_\mathcal{C}(gx)$ quickly,
we precomputed all projection chambers, that is, for each face $F$ we
computed an inequality representation of the polyhedral cone of all points
projected on $F$.  This precomputation is based on the following Lemma,
with the conversion to an inequality representation being done with the
software {\tt polymake} \cite{Polymake}.

\begin{lemma}
 \label{lem:proj-chambers}
 Let $K\subseteq [d]$ be of cardinality $k=|K|$.  Let
 \[
 F = \left\{ x\in\mathbb{R}^n \::\: (r_i,x) =0 \;\text{for all}\; i\in K,
 \;(r_i,x) \ge 0 \;\text{for all}\; i\in [d]\setminus K \right\}
 \]
 be the $k$-dimensional face of the fundamental chamber associated with
 $K$.  Then the set of all points $x\in V$ with projection
 $\pi_\mathcal{C}(x)\in F$ is the polyhedral cone spanned by $\{-r_i\::\:
 i\in K\}\cup \{s_i\::\: i\in [d]\setminus K\}$.
\end{lemma}

While our simulations offer strong evidence for the validity of
Conjecture~\ref{conj:main-lemma}, consideration of a single (or several)
random points unfortunately does not prove the conjecture.  The difficulty
lies in the fact that the orbits of different points in the fundamental
chamber behave differently with respect to projection chambers.  We refer
to Figure~\ref{fig:A3} where the orbits of two distinct points from the
fundamental chamber of $A_2$ are shown.  

\begin{remark}
 \label{rem:other-interpretations}
 The fact that the orbits behave differently with respect to projection
 chambers also makes it difficult to connect our interpretation of the
 coefficients with others that have appeared in the literature.
 For example, $|a_k|$ is also known to be the number of group elements in
 $\mathcal{W}$ that leave fixed all points of some linear space of
 dimension $n-k$; see \cite[Section 5]{ST}.  These two interpretations do
 not coincide, and it would thus be interesting to establish a bijection
 between the different partitions of $\mathcal{W}$.  The partition of
 $\mathcal{W}$ obtained through Conjecture~\ref{conj:main-lemma} will
 however depend on the orbit type of $x$.
\end{remark}

A computational approach towards a full proof of
Conjecture~\ref{conj:main-lemma} could proceed by decomposing the
fundamental chamber $\mathcal{C}$ of each of the six exceptional
irreducible reflection groups by orbit type.  For a given group
$\mathcal{W}$, define two points $x$ and $y$ equivalent if for all group
elements $g\in \mathcal{W}$, the projection $\pi_\mathcal{C}(gx)$ is in the
relative interior of the same face of $\mathcal{C}$ as
$\pi_\mathcal{C}(gy)$.  By a decomposition by orbit type we mean a
polyhedral subdivision of $\mathcal{C}$ that corresponds to equivalence
classes of the relation just defined.  A computer proof for $\mathcal{W}$
would then be complete upon consideration of the orbit of a single point
from each cone in the subdivision.

While such a computer proof may be feasible in a faster language than
Maple, it would be more desirable if a case-free proof of
Conjecture~\ref{conj:main-lemma} could be found, eliminating in
particular a separate treatment of the exceptional groups.

\bibliographystyle{amsalpha} 
\bibliography{reflection}

\end{document}